\documentclass[12pt]{amsart} 
\usepackage{amssymb}

\usepackage{euscript}
\let\script\EuScript
\let\cal\mathcal

\makeatletter \@mparswitchfalse \makeatother

\numberwithin{equation}{section}

\newtheorem{theorem}{Theorem}[section]
\newtheorem{lemma}[theorem]{Lemma}
\newtheorem{corollary}[theorem]{Corollary}

\theoremstyle{definition}
\newtheorem{example}[theorem]{Example}
\newtheorem{definition}[theorem]{Definition}

\def\rank{\operatorname{rank}}
\def\span{\operatorname{span}}
\def\range{\operatorname{Range}}
\def\HH{H\!H+H\!H}
\def\MM{M\!M}
\def\eqdef{\overset{\text{def}}=}
\def\O{\Omega}
\def\D{\mathbb D}

\def\N{\mathbb N}
\def\C{\mathbb C}
\def\CC{\script C}
\def\L{\script L}
\def\M{\script M}

\def\Z{\mathbb Z}
\def\R{\mathbb R}
\let\<\langle
\let\>\rangle
\let\phi\varphi
\let\epsilon\varepsilon

\let\\\cr
\let\Beta\beta
\def\eqref#1{(\ref{#1})}
\def\ps/{pseudomultiplier}
\begin{document}

\def\eqalign#1{\null\,\vcenter{\openup\jot
  \ialign{\strut\hfil$\displaystyle{##}$&$\displaystyle{{}##}$\hfil
      \crcr#1\crcr}}\,}
\def\eqalignbot#1{\null\,\vbox{\openup\jot
  \ialign{\strut\hfil$\displaystyle{##}$&$\displaystyle{{}##}$\hfil
      \crcr#1\crcr}}\,}
\def\eqaligntop#1{\null\,\vtop{\openup\jot
  \ialign{\strut\hfil$\displaystyle{##}$&$\displaystyle{{}##}$\hfil
      \crcr#1\crcr}}\,}

\date{May 23, 1996}

\title[Functions which are almost multipliers of Hilbert function spaces]
{Functions which are almost multipliers of\\Hilbert function spaces}
\author{J.~Agler}
\address{\hskip-\parindent J.~Agler\\
Department of Mathematics \\ 
University of California, San Diego \\
La Jolla CA 92093}

\author{N.~J.~Young}
\address{\hskip-\parindent N.~J.~Young\\
Department of Mathematics and Statistics\\
Lancaster University\\
Lancaster LA1 4YF \\
England}
\email{n.young@lancaster.ac.uk}

\thanks{Agler's research was supported by an NSF grant in Modern
Analysis.  Young wishes to thank
the University of California at San Diego and the Mathematical Sciences
Research
Institute for hospitality while this work was carried out.  Research at MSRI is
supported in part by NSF grant DMS-9022140. }

\begin{abstract}
We introduce a natural class of functions, the 
{\em pseudomultipliers}, associated with
a general Hilbert function space, prove an extension theorem
which justifies the definition, give numerous examples
and establish the nature of the $1$-pseudomultipliers of
Hilbert spaces of analytic functions under mild
hypotheses.
\end{abstract}

\maketitle

The function $1/z$ on the unit disc $\D$ is almost a
multiplier of the Hardy space $H^2$: it misses by only one dimension.
That is, there is a closed subspace of $H^2$ of codimension 1 
which is multiplied by $1/z$  into $H^2$. The same statement holds for
the
characteristic function of the point $0$. These are two key examples
of functions that we call \emph{pseudomultipliers} of Hilbert function
spaces. Now {\em multipliers} of the standard function spaces
have been much studied: it is a natural generalization to 
consider functions which fail to be multipliers by only finitely
many dimensions.  Moreover, for some familiar function spaces
one obtains in this way natural classes of functions.
For example, a well-known theorem of
Adamyan, Arov and Krein \cite{AAK} on $s$-numbers of Hankel operators
can be interpreted as a description of the \ps/s of $H^2$
(see Theorem \ref{taak} below).  They are the finite modifications
of functions of the form $f/p$ where $f$ is bounded and analytic 
in the unit disc and $p$ is a polynomial which does not vanish on the
unit circle.  The \ps/s of the Fock space
(see Theorem \ref{tfock} below) are the finite modifications of the
proper rational functions.  We address the question of what can
be said about 
the \ps/s of other popular function spaces.  

One can formulate the definition of a \ps/ of a function space
in very great generality.
In this paper we introduce the notion for the case of Hilbert 
function spaces.  This still covers a very wide variety of
spaces, but we are nevertheless able to make significant assertions
about them.  A major purpose of the paper is to give the
``correct" formulation of the notion of a \ps/.  To justify
our definition we establish an extension theorem (Theorem \ref{textn}).
This is the main result of the paper and shows that apparently
weaker variants of our definition give rise to essentially
the same objects.  Other goals are to illustrate the notion by means
of a range of examples and to show that a surprising amount
can be said about \ps/s of Hilbert spaces of {\em analytic} functions
under modest hypotheses.  Virtually all our results depend heavily on 
a key technical fact, Lemma \ref{l1.6}.
 
Our
investigation began as a study of certain interpolation problems.
It is known \cite{Ag1,Ag2,Q1} that the
classical result of Pick on interpolation by bounded analytic functions
\cite{P}
can be extended to certain other function spaces, some of them
having no connection with analyticity.
Could the same be true of the interpolation theorem
of Adamyan, Arov and Krein which generalises Pick's theorem?
This question was discussed but not resolved in \cite{Q2}.
To answer it we had to analyse the ``Pick kernel"
\begin{equation}
\label{pick}
(\lambda , \mu )\mapsto (1-\phi(\lambda )\bar{\phi}(\mu))k(\lambda ,\mu)
\end{equation}
corresponding to a reproducing kernel $k$ and a function $\phi$.
Pick's theorem concerns the positivity of this kernel on $\D$
when $k$ is the Szeg\H{o} kernel,
while the result of  Adamyan, Arov and Krein tells us what happens
when it has  $m$ negative squares on $\D$.  A study of the
same question for other kernels led us to the notion of a \ps/,
which proved to be a remarkably fruitful notion despite its great generality.
 We think of pseudomultipliers as being something like
meromorphic functions (plus point discontinuities), though they are
defined on arbitrary sets which need have no differentiable or even
topological structure - all information is contained in the kernel.
Regarding the original question, we shall show 
in a future paper that, for 
the most straightforward formulation,
the answer is negative: the Adamyan-Arov-Krein property characterizes
the Szeg\H{o} kernel.  However, \ps/s are interesting on their own account.
 
In this paper a \emph{kernel} on $\Omega$ will mean a function
$$
k:\Omega \times \Omega \to \C,
$$
where $\Omega$ is a set, satisfying
$$
k(\lambda,\mu)=k(\mu,\lambda)^- \quad\text{for all } \lambda, \mu\in
\Omega.
$$
We say that $k$ is \emph{nonsingular} if, for any $n\in \N$ and
$\lambda_1,\ldots, \lambda_n\in \Omega$, the $n\times n$ matrix
$[k(\lambda_i, \lambda_j)]^n_{i,j=1}$ is nonsingular.
We call $k$ a \emph{positive} kernel if every $[k(\lambda_i,
\lambda_j)]$ above is positive, and \emph{positive definite} if it is
both positive and nonsingular. For any non-negative integer $m$, we
say that a kernel $k$ has $m$ \emph{negative squares} if, for any
$n\in\N$ and $\lambda_1,\ldots,\lambda_n\in \Omega$, the $n\times n$
matrix $[k(\lambda_i, \lambda_j)]^n_{i,j=1}$ has no more than $m$
negative eigenvalues (counted with multiplicity), while for some
choice of $n, \lambda_1,\ldots, \lambda_n$, it has exactly $m$
negative eigenvalues.

If $k$ is a positive kernel on a set $\Omega$, there is
an associated Hilbert space $H(k)$ of functions on $\Omega$ for
which $k$ is the reproducing kernel \cite{Ar}. That is, for any 
$f \in H(k)$ and $\lambda \in \Omega$, we have
$$\< f, k_\lambda \>_{H(k)} =
f(\lambda),$$ 
where $k_\lambda \in H(k)$ is defined by
$k_{\lambda}(\mu) = k(\mu, \lambda)$ for $\mu \in \Omega$.
By a {\em Hilbert space of functions on a set $\O$} we shall mean a
Hilbert space $H$
whose elements are functions on $\O$ such that the point evaluations $h
\mapsto
h(\lambda)$ are continuous linear functionals on $H$ for all $\lambda
\in \O$.
Such a space has a reproducing kernel which is a positive kernel on
$\O$; it may
 fail to
be nonsingular.  If $H$ is a Hilbert space of functions on $\O$  
with kernel $k$ and $E$ is a subset of
$\O$ then the space of restrictions of elements of $H$ to $E$ is a
Hilbert space of
functions on $E$, whose kernel is the restriction of $k$ to $E$.  If $E$
is large enough
these two Hilbert spaces are naturally isomorphic, but for our purposes
they are to 
be regarded as two distinct Hilbert function spaces.  This point will be
important when
we discuss domains of definition of functions.

We say that a function $\psi:\Omega \to \C$ is a \emph{multiplier} of
$H(k)$ if $\psi f\in H(k)$ for every $f\in H(k)$. When $\psi$ is a
multiplier of $H(k)$, the mapping $f \mapsto\psi f$ is a linear
operator, which we denote by $M_\psi$. By the closed graph theorem,
$M_\psi$ is bounded. The \emph{multiplier norm} of $\psi$ is
$\|M_\psi\|$. The set of all multipliers of $H(k)$ is a subalgebra
of the algebra $\L(H(k))$ of all bounded linear operators on
$H(k)$. It is well-known \cite{BB, S} (and will follow from
Theorem \ref{textn} below) that, for any function $\psi: \Omega \to
\C$, the following two statements are equivalent:
\begin{enumerate}
\item[(a)]$\psi$ is a multiplier of $H(k)$ and $\|M_\psi\|\leq 1$, and
\item[(b)]the kernel $(\lambda, \mu) \mapsto
(1-\psi(\lambda)\bar{\psi}(\mu))k(\lambda, \mu)$ is positive on
$\Omega$.
\end{enumerate}

If $H$ is a Hilbert space and $u,v\in H$ we denote by $u\otimes v$ the
rank $1$ operator on $H$ given by $x\mapsto \< x,v\> u$. If $\CC$ is a
class of functions on a set $\Omega$ we say that $E\subset \Omega$ is
a \emph{set of uniqueness} for $\CC$ if $f,g\in\CC$ and $f|E=g|E$ imply
$f=g$ (here $f|E$ is the restriction of $f$ to $E$). If $k$ is a
positive definite kernel on $\Omega$ then $E\subset \Omega $ is a set
of uniqueness for $H(k)$ if and only if $\{k_\lambda: \lambda\in E\}$
spans a dense linear subspace of $H(k)$. For a space $H$ of functions
on $\Omega$ we shall have occasion to use the class
$$
\HH \eqdef \{f_1 f_2 + f_3 f_4: f_j\in H\hbox{ for } 1\leq j\leq 4\}.
$$
We denote by $s_j(T)$, for $j\geq 0$, the $s$-numbers of a bounded
linear
operator $T$ on a Hilbert space $H$:
$$
s_j(T) = \inf\{ \|T-R\|: R\in\L(H), \text{ rank }R\leq j\}.
$$
We remark that nonsingularity of a positive 
kernel $k$ on a set $\O$ can be expressed
in terms of its associated Hilbert space $H$ of functions.
It is easy to show that $k$ is nonsingular on $\O$ if and
only if, for every finite sequence $\lambda_1,\dots,\lambda_n$ of
points of $\O$, there is a function $h \in H$ such that $h(\lambda_1)
=1$ and $h(\lambda_j)=0$ for $2 \leq j \leq n$.  Another
equivalent condition is that $\{k_\lambda: \lambda \in \O\}$
is a linearly independent set in $H$.

\section{Pseudomultipliers of Hilbert Function Spaces}

The starting point for the notion of an $m$-pseudomultiplier
of a Hilbert space $H$ of functions on a set $\O$ is simply a function
$\psi$ on $\O$
such that, for some closed $m$-codimensional subspace $E$ of $H$, $\psi
E \subset H$.
However, as our archetype $1/z$ shows, we need to allow the possibility
that $\psi$ is
not defined on the whole of $\O$, and this gives rise to some delicate
issues.  
We
illustrate these issues by examples in order to motivate the formal
definition.
 In the two examples which
follow $k$ is the Szeg\H{o} kernel and $\Omega$ is $\D$, so that $H$
is~$H^2$.  We consider functions $\psi :D_{\psi} \subset \O \to \C$.
\begin{example}
\label{e1.2a}
$D_{\psi} = \D \setminus \{0\} \text{ and } \psi(\lambda)
=\lambda^{-m}$, where $m\in\N$. Here $z^mH^2$ has codimension $m$ in
$H^2$ and $
\psi$
multiplies  $z^mH^2$ into $H^2$.  Multiplication by $\psi$ determines a
bounded 
linear
operator of norm 1 from $z^mH^2$ into $H^2$; we say that $\psi$ is an 
$m$-\ps/ of $H^2$ and
that $s_m(\psi)=1.$
\end{example}

\begin{example}
\label{e1.3a}
$D_\psi=\D$ and $\psi$ is given by $\psi(0)=1$ and
$\psi(\lambda)=0$ if $\lambda\neq 0$.
Here multiplication by $\psi$ determines the zero operator from the
1-codimensional
subspace $zH^2$ of $H^2$ into  $H^2$.  Thus $\psi$ is a 1-\ps/ and
$s_1(\psi)=0$
.
However $\psi$ is not a $0$-pseudomultiplier.
\end{example}

What should we assume about the domain of an $m$-\ps/?  It should be in
some sense
``almost all" of $\O$, but in what sense?  There are various choices one
could make; an
important purpose of this paper is to show that all the natural
assumptions about the
domain lead to essentially the  same objects, and these are functions
which are 
defined
everywhere on $\O$ with the exception of at most $m$ points. 
Accordingly we shall
incorporate this requirement into the definition, and in the following
section we shall
justify our decision by showing that weaker assumptions lead to the
same class of
functions.

A $1$-\ps/ of $H^2$, then, will be a function defined on the complement
of at most one 
point in $\D$.  Consider the function $\psi$ which is identically zero
on
$\D\setminus \{0\}$.  This $\psi$ satisfies $\psi zH^2 = \{0\} \subset
H^2.$  Should we
regard $\psi$ as a \ps/?  And as distinct from the function which is
identically 
zero on
all of $\D$?  We prefer to regard these two functions as the same, and
so will include
in our definition a condition to ensure that if a \ps/ is undefined at a
point of $\O$
then it is in a precise sense not definable at that point.

Now consider a fixed Hilbert space $H$ of functions on a set $\O$ with
nonsingular
kernel $k$.  For $D \subset \O$ and any function $f$ on $\O$ we denote by
$f|D$ the restriction of $f$ to $D$, and we denote by $H|D$ the set 
$\{h|D: h \in H\}$.
For $\psi : D_\psi\subset \O \to \C$ and $h \in H$ we denote
by $\psi h$ the pointwise product of $\psi$ and $h|D_\psi$, so that
$\psi h$ is a
function on $D_\psi$.  For $E \subset H$ we denote by $\psi E$
the set $\{\psi h: h \in E\}$, so that $\psi E$ is a set of functions on
$D_\psi$.
 If $E$ is a closed linear subspace
of $H$, if $D_\psi$ is a set of uniqueness and
 $\psi E \subset H|D_\psi$ we define the linear 
operator $M_{\psi,E} : E \to H$
by
$M_{\psi,E}h = g$ where $g$ is the unique element of $H$ which extends
$\psi h$.
  By the
closed graph theorem $M_{\psi,E}$ is continuous.  
If $D_\psi$ is a set of uniqueness of $H$ 
we shall say that $\alpha
\in \O$ is a {\em
singularity} of $\psi$ (relative to $H$)
if there exists $h \in H$ such that $h(\alpha) = 0$ and
$\psi h$ has an extension to a function $g \in H$ (necessarily unique)
such that $g(\alpha) \neq
0$.  Of course this can only happen if $\alpha \not \in D_\psi$.

\begin{definition}
\label{d1.1}
Let $\psi$ be a complex-valued function on a subset $D_\psi$ of
$\Omega$ and let $m\in\Z^{+}$, the set of non-negative integers. We
say that $\psi$ is an \emph{m-pseudomultiplier} of $H$ on $\Omega$ if
\begin{enumerate}
\item $D_\psi$ is a set of uniqueness for $H$;
\item $\O \setminus D_\psi$ consists of singularities of $\psi$ 
and contains at most $m$ points;
\item there is a closed subspace $E$ of codimension $m$ in $H$ 
such that  $\psi E \subset H|D_\psi$.
\end{enumerate}

For an $m$-\ps/ $\psi$ of $H$ we define
$s_m(\psi)$ to be the
infimum of $\| M_{\psi,E}\|$ over all closed $m$-codimensional subspaces
$E$ of 
$H$ such
that $\psi E \subset H|D_\psi$.
We say that $\psi$ is a {\em \ps/} of $H$ if it is an $m$-\ps/ of $H$ for
some non-negative integer $m$.
\end{definition}
\noindent {\bf Notes} 1. In the case of a finite set $\O$ 
no proper subset is a set of
uniqueness (recall
that we have supposed the kernel $k$ nonsingular), and so \ps/s are
everywhere defined
on $\O$.  Every $\psi: \O \to \C$ is an $m$-\ps/ and $s_m(\psi)$ is the
$m$th singular
value of the multiplication operator $M_\psi$.

2. A $0$-\ps/ $\psi$ is the same as a multiplier of $H$, and $s_0(\psi) =
s_0(M_{\psi}).$ 

3.  If $m \ge 1$, $\psi$ is defined at all but $m-1$ points of
$\Omega$ 
and $\psi$ is an $m$-\ps/ but not an $(m-1)$-\ps/ then there is a {\em
unique} $E$ with the stated properties, and so $s_m(\psi) = \|
M_{\psi,E}\|$.  For if
$E_1,~E_2$ are distinct and both have all the stated 
properties then $E_1 +E_2$ is a closed subspace of $H$ of
codimension at most $m-1$ and $\psi(E_1 +E_2) \subset H|D_\psi$.  
That is, $\psi$ is an $(m-1)$-\ps/.

In particular, if a $1$-\ps/ $\psi$ is defined on the whole of $\O$
then unless it is a multiplier, there is a unique $1$-codimensional
subspace $E$ of $H$ such that $\psi E \subset H$.
 However,  it can happen for a 1-\ps/ $\psi$ for which $D_\psi \neq \O$
that $\psi H \subset H|D_\psi$ and yet $\psi$ cannot be extended
to a multiplier of $H$: see Example \ref{e1.6} below.

4. A related definition of $s_m(\psi)$ and an analysis in the case of
``Blaschke kernels'' $k$ is in \cite{Q2}.

5. For $\psi$ as in the Definition, $\O \setminus D_\psi$ consists
precisely of the singularities of $\psi$.

The notion of a \ps/ has as much generality as that of a Hilbert
function space, and so it is not to be expected that we can give
a very precise taxonomy of \ps/s in general.  We do want to give
a wide range of examples to illustrate the notion; to introduce a 
little order into these examples we shall take a preliminary
step towards a classification of $1$-\ps/s.  For this purpose we
shall make an additional (but mild) assumption on the \ps/ $\psi:
D_\psi \to \C$, to wit that $D_\psi$ is a set of uniqueness
for the class $\HH$ of functions on $\O$.  A simple distinction
is between the $\psi$ for which  $D_\psi=\O$ and those for
which $\O \setminus D_\psi$ contains exactly one point.  We
shall say that $\psi$ has {\em defect} 0 or 1 according as the first
or second of these holds.  Consider first the case of defect 1.
There exists a unique $\alpha \in \O \setminus D_\psi$, and there
are $u, g_1 \in H$ such that $u(\alpha) = 0,~\psi u =g_1|D_\psi,$
and $g_1(\alpha) \neq 0$.  There is also a closed $1$-codimensional
subspace $E$ of $H$ such that $\psi E \subset H|D_\psi$.  We claim that
$$
E = k_\alpha^\perp =\{h \in H: h(\alpha)=0\}.
$$
For suppose that $h \in H$ is such that $\psi h \in H|D_\psi$ -- say
$\psi h = g_2 |D_\psi,~ g_2 \in H$.  On $D_\psi$ we have
$$
g_1h - ug_2 = \psi uh -\psi uh =0,
$$
and so, by the hypothesis on $\HH$, $g_1h = ug_2$ on $\O$.  Since
 $u(\alpha) = 0$ and $g_1(\alpha) \neq 0$ we have $ h(\alpha)=0$.
Thus $E \subset  k_\alpha^\perp$.  Since $E$ has codimension 1
it must be that $E =  k_\alpha^\perp$.  We have shown that,
subject to the assumption on $\HH$, if $\psi$ is a $1$-\ps/ of
defect 1 and $D_\psi = \O \setminus \{\alpha \}$, then
$\psi  k_\alpha^\perp \subset H|D_\psi$, and moreover
$E =k_\alpha^\perp$ is the only closed 1-codimensional
subspace of $H$ such that $\psi E \subset H|D_\psi.$

Now suppose that $\psi$ is of defect 0, i.e. $\psi$  is everywhere defined
on $\O$.
We identify three subcases.
\begin{enumerate}
\item $\psi$ is a multiplier of $H$.
\item  $\psi$ is not a multiplier and there is a point $\alpha \in \O$
such that
$\psi k_\alpha^\perp \subset k_\alpha^\perp.$
\item  $\psi$ is not a multiplier and there is a closed subspace $E$ of
codimension $1$
in $H$ such that $\psi E \subset H|D_\psi$ and $E$ has no zeros in $\O$.
\end{enumerate}

These three cases do exhaust the 1-\ps/s of defect 0.  For if $\psi$ is
not a multiplier
then, by Note 3 above, there is a unique closed 1-codimensional subspace
$E$ of 
$H$ such
that $\psi E \subset H|D_\psi$.  If $E$ is of the form $k_\alpha^\perp$ for some
$\alpha \in
 \O$
then we must have  $\psi k_\alpha^\perp \subset k_\alpha^\perp$ since
$\psi(\alpha)$ is defined,
and we are in case 2.  The remaining possibility is that $E$ is not
$k_\alpha^\perp$ 
for any
$\alpha$, and case 3 applies.

Let us exhibit examples of all these cases.  Defect 1 is exemplified by
Example
\ref{e1.2a} ($\psi(z) = 1/z$).  Defect 0, case 1 is multipliers --
constant functions will
do.  Example \ref{e1.3a} (the characteristic function of a point)
illustrates defect 0,
case 2.  Defect 0, case 3 occurs in the following.
\begin{example}
\label{ewrongO}
Let $k$ be the Szeg\H{o} kernel
restricted to the disc $\Omega = \{\lambda \in \C: |\lambda -
\frac{1}{2}| < \frac{1}{2}\}$, take $D_\phi = \Omega $ and
$\phi(\lambda) = \lambda^{-1}$. Then $s_1(\phi) = 1$.  
As a Hilbert space $H$ is the usual Hardy space $H^2$, but as a function
space
it is a space of
functions on $\Omega$.   Multiplication by $\phi$ acts
contractively from $\boldsymbol1^\perp$ into $H$, but the elements of
$\boldsymbol1^\perp$ have no common zero in~$\Omega $.
\end{example}

The latter example looks somewhat unnatural, because we took $k$ on the
``wrong'' domain $\Omega $, but we are working in such generality that
such displeasing examples are included. In the next section we
shall make further assumptions on $k$ in order to exclude such
phenomena, and we shall be able to make stronger assertions about the
structure of pseudomultipliers.
Let us return to defect 0, case 2: $\psi k_\alpha^\perp \subset k_\alpha^\perp$. 
Note that 
$\psi $
is defined at $\alpha$, but that  the value of $s_1(\psi)$ is
independent of
$\psi(\alpha)$.  We could redefine $\psi(\alpha)$ without affecting the
fact that $\psi$
is a 1-\ps/ or the value of   $s_1(\psi)$.  We shall call such a
phenomenon an {\em
anomaly} of $\psi$.  Our example \ref{e1.3a}, the characteristic
function of the point
$0$, has an anomaly at $0$.  In this example it is possible to redefine
$\psi(0)$ to
obtain a multiplier of $H$.  This is one easy way to generate 1-\ps/s of
defect 0, case 2:
we will always obtain either such an object or a multiplier if we take a
multiplier and change its value at a single point.  However, not all
instances of defect
0, case 2 arise in this way.
\begin{example}
\label{eanom}
Let $\O = \D$ and let $H = zH^2 \oplus \C g$ where $g(z) = (z-1)^{-1}$
and the norm in $\C g$ is given by $\| \lambda g\| = | \lambda |.$
  Let $\psi(z) =
z$ if $ z \neq 0$ and let $\psi(0) =c$.  Then $\psi H(0) \subset H(0)$,
so that 
$0$ is
an anomaly of $\psi$, but there is no choice of $c$ which makes $\psi$ a
multiplier of
$H$.   For suppose that $\psi g \in H$.  Then
$$
\psi g = zf + \lambda g \text{  for some } f \in H^2, \quad \lambda \in
\C.$$
On putting $z=0$ we find that $\lambda = c$, and so$$
f(z) = \frac{\psi(z) - c}{z}g(z) = \frac{z-c}{z(z-1)} \text{ for } z
\neq 0.$$
This contradicts $f \in H^2$, whatever the choice of $c \in \C$.
\end{example}

Let us extend the terminology introduced above.  For a closed subspace
$F$ of $H$ we
denote by $P_F$ the orthogonal projection operator from $H$ to $F$.
\begin{definition}
\label{d1.2}
Let $\psi$ be a $m$-\ps/ of $H$.
We say that $\psi$ is of {\em defect} $j$ if the complement of the 
domain of $\psi$ contains $j$ points.  We say that $\alpha \in \O$ is
an {\em anomalous point} or an {\em anomaly} of $\psi$ if $\psi
k_\alpha^\perp \subset  k_\alpha^\perp$.
\end{definition}

Here are further examples of the variety of behaviour of \ps/s.
\begin{example}
\label{e1.2b}
A pole which is not a singularity.
 $H =H^2$, $D_\phi = \D$ and
$$
\phi(\lambda)=
\begin{cases}
\lambda^{-m}&\hbox{if }\lambda\ne0,\\
c&\hbox{if }\lambda=0,
\end{cases}
$$
where $m\in \N$ and $c\in \C$.
This example differs from Example \ref{e1.2a} in that $\phi$ is defined
at $0$.
Here $\phi$ multiplies $z^{m+1}H$, but no $m$-codimensional subspace of
$H$, into $H$.
It is thus an $(m+1)$-\ps/ but not an $m$-\ps/.  It has defect $0$.
\end{example}
In connection with examples \ref{e1.2a} and \ref{e1.2b} we pose the
question: is
 there a
natural notion of the multiplicity of a singularity of a \ps/ for a
general kernel?

\begin{example}
\label{ell2}
The simplest of all infinite-dimensional Hilbert
function spaces is $\ell^2$ as a space of functions on $\N$ in the
natural
way.  Its kernel is $k(\lambda, \mu)= \delta_{\lambda \mu}$, the
Kronecker
symbol.  No proper subset of $\N$ is a set of uniqueness for $\ell^2$,
and
so \ps/s are defined on all of $\N$.  In fact the \ps/~algebra of
$\ell^2$ is $\ell^{\infty}$, which is also the multiplier algebra.  For
suppose
there is an unbounded $m$-\ps/~$\psi$.  We may pick an increasing
sequence
($n_j$) in $\N$ such that $| \psi (n_j)| \geq j$ and a space $E$ of
codimension $m$ in $\ell^2$ such that $\psi E \subset \ell^2$.  
Let $P_N$ be the
orthogonal projection operator from $\ell^2$ to span$\{k_1, k_2, \dots,
k_N\}$.  Choose $N$ so large that $P_N E^\perp$ has dimension $m$. 
Then,
if $f_1, \dots, f_m$ is a basis of $E^{\perp}$, the vectors $P_Nf_1,
\dots,
P_Nf_m$ are linearly independent.  For $n>N$ define $x(n)$ to be $1/j$
if
$n=n_j,$ 0 otherwise.  Pick $\xi \in$ span$\{k_1, \dots, k_N\}$
such that 

$$
<\xi, P_Nf_j> = - \sum_{n=N+1}^{\infty} x(n) \bar{f}_{j} (n),
\text{ for } 1 \le j \le m.
$$

This is possible on account of the linear independence of the $P_Nf_j$. 
Now set $x(n) = \xi(n)$ for $1 \le n \le N$: then $x{\perp}f_j$ for 
$1 \le j \le m$. Thus
$x \in  E$ and hence $\psi x \in \ell^2$.  However
$| \psi x (n_j)| \geq 1$ for infinitely many points $n_j$, a
contradiction.
Thus $m$-\ps/s of $\ell^2$ are bounded.  
\end{example}

\begin{example}[The Fock space $\Phi$]
\label{e1.3b}
This is the space with kernel $k(\lambda,\mu)=e^{\lambda\bar{\mu}}$ on
$\C$. It is also known as the Fischer space and the Bargmann--Fock
space. It comprises all entire functions $f$ such that
$$
\int_{\C} |f(z)|^2 e^{-|z|^2} \,dA < \infty
$$
where $dA$ is area measure. Alternatively it can be
described as the space of functions $f(z)=\sum^\infty_{n=0}a_n z^n$
such that $\sum^\infty_{n=0} |a_n|^2 n! < \infty$. An orthonormal
basis of $\Phi$ is given by the functions
$\bigl(z^n / \sqrt{n!}\,\bigr)^\infty_{n=0}$. The interest of $\Phi$
here is
that it admits no nonconstant multipliers, but it has plenty of
nontrivial pseudomultipliers. 
Indeed, a non-trivial Hilbert space of entire functions
cannot have a non-constant multiplier.  For suppose that
$\phi$ is such a multiplier.  It is easily shown that $\phi$ is entire,
and so, by Liouville's theorem, $\phi(\C)$ is unbounded.  Since
$H$ is non-trivial its kernel $k$ satisfies $k_\lambda \neq 0$ for
$\lambda$ in a dense subset of $\C$.  The relation $M_\phi^* k_\lambda=
\bar{\phi}(\lambda)k_\lambda$ shows that $\phi(\C) \subset \sigma(M_\phi)$,
contradicting boundedness of $M_\phi$.
In particular therefore the only multipliers of the Fock space 
are the constants,
but $\phi:
\C\setminus\{0\}\to\C$ with $\phi(\lambda)=\lambda^{-1}$ is a
$1$-pseudomultiplier: it multiplies the $1$-codimensional subspace 
$\{h \in H: h(0)=0 \}$ of $H$ into $H$.  We shall discuss the \ps/s of $\Phi$
further in Sec. 3 below.
\end{example}
$H^2$ has the property that if $\psi$ is a multiplier then $s_m(\psi) =
s_0(\psi)$ for all $m \in \Z^+$.  A simple example of a space
which does not have this property is $\ell^2$, as a space of functions on
$\N$.  Here $\psi(n)= 1/n$ defines a multiplier of $\ell^2$; clearly
$s_m(\psi) = 1/(m+1)$.  

One might imagine that a meromorphic function defining a 1-\ps/ would
have a single, simple pole.  This is true for nonsingular kernels,
as will follow from
Theorem \ref{mplcty} below, but it is not true if we relax for a moment
our standing assumption that $H$ has nonsingular kernel, as
the following example shows.  
\begin{example}
\label{e1.7}
Let $m \geq 2$ and let $H$ be the  closure in $H^2$ of the space of polynomials in
$\lambda^m$, as a space of
functions on $\D$, and let $\phi(\lambda)= \lambda^{-m}$ on $\D \setminus \{0\}.$  Then
$\phi$ multiplies $\boldsymbol{1}^\perp$ into $H$, 
hence is a $1$-pseudomultiplier, but has a pole of multiplicity $m$ at $0$.
\end{example}

Examples \ref{e1.2a} and \ref{e1.3b} suggest a conjecture: if $H$ is a
Hilbert space of analytic functions on a domain $\O$ with nonsingular
kernel, then for any multiplier $\psi$ of $H$ and any $\alpha \in \O$
then the function  $\lambda \mapsto \psi(\lambda)/(\lambda - \alpha)$
 is a $1$-pseudomultiplier of $H$.
 However, this is not so, even if $H$ is $z$-invariant.
\begin{example}
\label{eHRS}
Let $\Omega = \D$ and $D_\phi= \D \setminus \{0\},
~\phi(\lambda) = \lambda^{-1}$.
Let $H$ be a $z$-invariant subspace of the 
Bergman space $L^2_a(\D)$ with nonsingular kernel such that
 $zH$ has codimension 2 in $H$.  Here $z$ denotes the
operation of multiplication by the independent variable.  It is shown in \cite{H}
how  such a subspace can be constructed.  
One can take $H$ to be of the form $I(A) + I(B)$
where $A,~B$ are certain zero sequences of  $L^2_a(\D)$
and the closed subspace  $I(A)$ is
$\{f \in L^2_a(\D): f = 0 \text{ on } A\}$.  
It is shown that $A,~B$ can be chosen so
that $A \cup B$ is a ``sampling sequence" for  $L^2_a(\D)$, which ensures that 
$I(A) + I(B)$ is closed. If $\lambda_1,\dots,\lambda_n$ are distinct points of
$\D$ one can find $f \in  I(A) + I(B)$ such that $f(\lambda_1) \neq 0$; then
the function
$\prod_{j=2}^n (\lambda -\lambda_j)f(\lambda) \in H$ vanishes at 
$\lambda_2,\dots,\lambda_n$ but not at $\lambda_1$ and so $H$ has nonsingular
kernel.
Now suppose that $\phi$ is a
$1$-pseudomultiplier of $H$.  By Definition \ref{d1.1} there is a
subspace $E$ of codimension one in $H$ and a bounded  linear operator
$L:E \to H$
such that $(Lh)(\lambda) = \lambda^{-1}h(\lambda)$ 
for all $h \in E$ and all but one $\lambda \in
\D \setminus \{0\}$.  Thus each $h \in E$ agrees with an element of $zH$ except for at
most two points, and so belongs to $zH$.  Thus $E \subset zH$, 
a contradiction since $E$
has codimension one and $zH$ has codimension two.
Thus $\phi$ is not a $1$-pseudomultiplier.
\end{example}
By a slight modification of the above construction one can 
find a $z$-invariant $H \subset L^2_a(\D)$
with nonsingular kernel
such that $\phi(\lambda) = \lambda^{-1}$ is not an $m$-pseudomultiplier for any integer
$m$.  It is shown in \cite{HRS} that $H$ can be chosen so that $zH$ has infinite
codimension in $H$.  It follows that $\phi$ is not a pseudomultiplier.

Most of our natural-looking examples of \ps/s relate to
Hilbert spaces of analytic functions, so it is of interest to 
show that other types of space can also have non-trivial
\ps/s.

\begin{example}
\label{esobolev}
A Sobolev space.  Let $W$ be the space $\{f \in L^2(0,1):
f^\prime \in L^2(0,1)\}$, with its usual inner product,
as a space of functions on $[0,1]$.  The characteristic
function of a point is a \ps/ of this space, but a much more
interesting example is the function $\phi(t)= \sqrt{t}$ .
Since $W$ is an algebra with identity 
and $\phi \not \in W$ it follows
that $\phi$ is not a multiplier of $W$, but $\phi$ {\em is}
a $1$-\ps/ of $W$.  In fact $\phi k_0^\perp \subset W$.
Consider $f \in W$ such that $f(0)=0$.  Then for $t \in (0,1)$
$$
|f(t)| = \left| \int_0^t f^\prime(x) dx \right| \leq \sqrt{t} ||f||_W,
$$
so that  $\phi^\prime f$ is bounded on $(0,1)$.  It follows
that $\phi f \in W$.  Hence $\phi k_0^\perp \subset W$ as
claimed.
\end{example}
We shall need the following simple fact.
\begin{theorem}
\label{mltplr}
If $\phi$ is a \ps/ with domain $D_\phi$ and $\phi$ is the restriction
of a multiplier of $H$ then $D_\phi = \O$ and $\phi$ is a multiplier of $H$.
\end{theorem}
\begin{proof}
We have to show that $\phi$ cannot have any singularities 
in $\O$.  Suppose, on the contrary, that $\phi$ 
 is the restriction
of a multiplier $\theta$ and
 that there exist $\alpha \in \O$ and 
$h \in H$ such that $h(\alpha) = 0$ and $\phi h$ has an extension
$g \in H$ such that $g(\alpha) \neq 0$.  Now $\theta h \in H$ extends
$\phi h$, and so $\theta h = g$.Thus
$$
0 = \theta(\alpha)h(\alpha) = g(\alpha) \neq 0,
$$
a contradiction.
\end{proof}

\section{Extension of \ps/s and the Pick kernel}

In our definition of an $m$-\ps/~$\phi$ of a Hilbert space $H$ of functions
on a set $\Omega$ the main idea is that $\phi$ multiply a closed
$m$-codimensional subspace of $H$ into $H$.  This is natural enough,
as is the requirement that the domain $D_\phi$ of $\phi$ be a set
of uniqueness for $H$; otherwise $\phi$ is more appropriately studied
via its action on the proper closed subspace of $H$ spanned by $k_\lambda,~
\lambda \in D_\phi$.  However, we have also imposed the condition 2 of
Definition \ref{d1.1}, that $\O \setminus D_\phi$ consist of 
singularities of $\phi$ and have at most $m$ points.  This is less
obviously appropriate.  Why make such a strong assumption on $D_\phi$?
It would make perfect sense simply to delete condition 2.  Let us 
temporarily call functions satisfying the resulting definition
{\em weak m-\ps/s}.  On the face of it these will comprise a larger
class than the $m$-\ps/s.  It is a remarkable fact that under a
modest hypothesis on $D_\phi$ every weak $m$-\ps/ can be extended with
preservation of $s_m(\cdot)$ to a function which is defined 
on all but at most $m$ points of $\O$ and which is an $m$-\ps/.

Observe that if $D_{\phi}$
is a set of uniqueness for $H$ then $H|D_{\phi} \eqdef \{h|D_{\phi}:
h \in H\}$ is a space of functions on $D_{\phi},$ and $h
\mapsto h|D_{\phi}$ is a Hilbert space isomorphism of $H$ and
$H|D_{\phi}$.  To say that $\phi$ is a weak $m$-\ps/ of $H$ 
is the same as saying that $\phi$ is an $m$-\ps/ of  $H|D_{\phi}$.
Our main aim in this section is to show that an $m$-\ps/ of  $H|D_{\phi}$
extends  with
preservation of $s_m(\cdot)$
to an $m$-\ps/ of  $H$.
 In fact this is not true in full generality, even
for multipliers $(m=0)$, as we show in Example \ref{e1.6} below, but the 
extra hypothesis we need to obtain such an extension holds for nearly all
spaces of interest.

\begin{theorem}
\label{textn}
 Let $H$ be a Hilbert space of functions on a set $\Omega$ with
nonsingular kernel $k$, 
let $D_{\phi} \subset \Omega$ be a set of uniqueness
for
$H$ and for $HH+HH$.  The following are equivalent for a function
$\phi: D_{\phi} \rightarrow \C$.
\begin{enumerate}
\item $\phi$ is an $m$-\ps/ of $H|D_{\phi}$;
\item There exists $t \geq 0$ such that the kernel
\begin{equation}
\label{pickt}
(\lambda, \mu) \mapsto \left(t^2 - \phi (\lambda) \bar{\phi} (\mu)
\right) k(\lambda,\mu)
\end{equation}
 has at most $m$ negative squares on $D_{\phi}$;
\item There is an $m$-\ps/ $\psi$ of $H$ which extends $\phi$.
\end{enumerate}
Moreover, when these conditions hold, 
$$
s_m(\phi) = s_m(\psi) = \inf t
$$
over all $t$ for which \eqref{pickt} has at most $m$ negative squares.
\end{theorem}

Before we prove the theorem let us look at an example which shows how
the conclusions
can fail if the hypothesis about $\HH$ is omitted.

\begin{example}
\label{e1.6}
Let $\Omega=\N$, let $u=\sqrt{3}(2^{-n})^\infty_{n=1}$, so that
$u$ is a unit vector in $\ell^2$, and let $H$ be the orthogonal
complement of $u$ in $\ell^2$. The reproducing kernel of $\ell^2$ is
$(\lambda,\mu)\mapsto\delta_{\lambda \mu}$  and
so the reproducing kernel in $H$ is given by orthogonal projection
onto $H$:
$$
k_\lambda(\mu) = [\delta_\lambda  - \< \delta_\lambda, u \> u] (\mu) =
                 \delta_{\lambda \mu} - \frac{3}{2^{\lambda + \mu}}.
$$

This example is simple, but it has some interesting properties.  It
illustrates pathologies which can occur when \ps/s are studied in
excessive generality.  

(a) Only the constant sequences are multipliers of $H$, but there is a
rich class of \ps/s.  In fact $\psi$ is a \ps/ of $H$
$\Leftrightarrow \psi$ is a $1$- \ps/ of $H$
$\Leftrightarrow$ either $\psi \in \ell^{\infty}$ or $\psi$
is a bounded non-constant function on $\N\setminus \{\alpha\}$ 
for some $\alpha \in \N$.

(b) For any $\alpha \in  \N$ and any non-constant bounded
function $\phi: \N\setminus\{\alpha\} \rightarrow \C$, $\phi$ is a
multiplier of
$H|\N\setminus\{\alpha\}$ but $\phi$ has no extension to a multiplier of
$H$.  

(c) For $\phi$ as in (b) and $t=\sup \{| \phi _n|: n \in
\N\setminus\{\alpha\}\}$, the Pick kernel
$$
(\lambda, \mu) \mapsto \left(t^2- \phi (\lambda)
  \bar{\phi}(\mu)\right) k(\lambda,\mu)
$$
 is positive but $\phi$ does not extend to a multiplier of $H$.  

(d) There are functions $\psi$ such that $\psi H \subset H|D_\psi$ but $\psi$ is
not a multiplier of $H$, nor even extendable to a multiplier of $H$.

{\bf Notes.}  Since $\N\setminus\{\alpha\}$ is a set of uniqueness for $H$,
for any
$\alpha \in  \N$, (b) shows that the implication $1
\Rightarrow 3$ of Theorem \ref{textn} fails in general, even when $m=0$.
 Likewise (c) shows that $2 \Rightarrow 3$ fails.  Here $\HH =
\ell^{1}_{\N}$, so that $\N\setminus\{\alpha\}$ is not a set of
uniqueness for $\HH$.

Let $\phi=(\phi_n)^\infty_{n=1}$ be a multiplier of $H$. Then $\phi x
\in H$, where 
$$x= (1,0, \ldots, 0, -2^{n-1}, 0, \ldots),$$ the
$-2^{n-1}$ being in the $n$th place. It follows that $\phi_n= \phi_1$.
Hence all multipliers of $H$ are constant.

Now consider $\psi = (\psi_n)^{\infty}_{n=1} \in$
$\ell^{\infty}_\N$. It is routine to show that if $f$ is the orthogonal
projection of $\psi u$ onto $u^\perp$  in $\ell^2$ then $\psi f^\perp
\subset u^\perp$.  
Hence $\psi$ is a $1-$\ps/ of $H$.  Thus the set of 1-\ps/s
of $H$ contains $\ell^{\infty}$; it also contains many other functions. 
Let $\alpha \in \N$ and let $D=\N\setminus\{\alpha\}$.  Then $D$ is a
set of uniqueness for $H$, since if $x \in H$ and we know $x_n$
for $n \neq \alpha$ we can recover $x_{\alpha}$ from the relation
$x \perp u$.  Let $\phi: D \rightarrow \C$ be bounded and
non-constant. 
Then
$\phi$ is a $1$-\ps/ of $H$.  Indeed, $\phi k_\alpha^\perp \subset H$, since
$\phi x \in \ell^2_D$ for any $x \in  H$, and every
element of $\ell^2_D$ has a unique extension to an element of $H$. 
Moreover, $\phi$ has a singularity at $\alpha$ in the sense of
Definition $\ref{d1.1}$.  Since $\phi$ is non-constant one can construct
$h \in H$  such that $h(\alpha)=0$ but $(\phi h)(\alpha)
\neq 0$.  For example, if $\alpha \not \in \{1,2\}$ and $\phi(1) \neq
\phi(2)$
we may take 
$h(1) = \phi(2)$, $h(2) = -2 \phi(1)$, $h(n)=0 \text{ for } n > 2$.
Hence $\phi$ is a $1$-\ps/ of $H$.

We have shown that the $1$-\ps/s of $H$ contain $\ell^{\infty}$ and all
non-constant bounded functions on complements of singleton subsets of
$\N$.  In fact these are all the \ps/s of $H$.  For let $\psi$ be
an $m$-\ps/.  Suppose $\psi$ is defined everywhere on $\N$ except at
a singularity
$\alpha$ (the domain of $\psi$ has to be a set of uniqueness, and so can
omit at most one point).  Pick a closed $m$-codimensional subspace $E$
of $H$ such that $\psi E \subset H$.  Then
$E$ has finite codimension as a subspace of $\ell^2$ and satisfies $\psi
E \subset \ell^2$.  It follows that $\psi$ is bounded (see Example
\ref{ell2}).  
It is easy to see that if $\psi$ is constant it cannot have
a singularity at $\alpha$.  Thus $\psi$ is a bounded non-constant
function on $\N\setminus\{\alpha\}$.  A similar argument shows that if
$\psi$ is
everywhere defined on $\N$ then $\psi \in \ell^{\infty}$.  This
completes the proof of (a).

(b) Let $D=\N\setminus\{\alpha\}$ for some $\alpha \in \N$.  We have
$H|D$ $=$ $\ell^2_D$, so that the multipliers of $H|D$ are just the
elements of $\ell^{\infty}_{D}$.  However, only the constant sequences
extend to multipliers of $H$.

(c) It will be seen below that the implication $1 \Rightarrow 2$ of
Theorem $\ref{textn}$ does not depend on the domain of $\phi$ being a set
of uniqueness for $\HH$, and so $\phi$ has the property described.

(d) Let $\psi$
be a non-constant bounded function on $\N\setminus\{\alpha\}$, where
$\alpha
\in \N$.  For any $h \in H$, $\psi h$ is an $\ell^2$
sequence indexed by $\N\setminus\{\alpha\}$, and it has a (unique)
extension to
an element of $H$.  Thus $\psi H \subset H|D_\psi$.  However, $\psi$ is not a
multiplier of $H$ since it is not defined at $\alpha$.  It does not even
extend to a multiplier of $H$ since it is non-constant.

\end{example}

{\it Proof of Theorem} \ref{textn} (3) $\Rightarrow$ (1) is immediate
from the
fact that $h \mapsto h|D_{\phi}$ is an isomorphism of Hilbert
spaces.  

(1) $\Rightarrow$ (2).  Let $E$ be a closed $m$-codimensional
subspace of $H$ such that $\phi E  \subset H|D_\phi$, 
so that $M \eqdef M_{\phi,
E}$ is a bounded linear operator from $E$ to $H$.
 
For $\lambda\in D_\phi$ and $h\in E$ we have
$$
\eqalign{
\<M^*k_\lambda ,h\>_{E} &= \<k_\lambda , Mh\>_H = (Mh)(\lambda)^-
                              = \bar{\phi}(\lambda)\bar{h}(\lambda)\cr
                              &= \bar{\phi}(\lambda)\< k_\lambda,h\>_H
                              =
\<\bar{\phi}(\lambda)P_{E}k_\lambda, h\>_{E},}
$$
where $P_{E}:H\to E$ is orthogonal projection. That is,
$$
M^*k_\lambda = \bar{\phi}(\lambda)P_{E}k_\lambda 
\quad\hbox{for } \lambda\in D_\phi.
$$
For any $t\in\R$ and $\lambda,\mu\in D_\phi$ we have
$$ \eqalign
{
\<(t^2 - MM^*)k_\mu,\,k_\lambda\> &= t^2k(\lambda,\mu) -
\<M^*k_\mu,\,M^*k_\lambda\> \cr  
                                &= t^2k(\lambda,\mu) -
\<\bar{\phi}(\mu)P_{E} k_\mu ,\, \bar{\phi}(\lambda)P_{E}
k_\lambda\> \cr                               
&=t^2k(\lambda,\mu)-\phi(\lambda)\bar{\phi}(\mu)
 \< P_{E}k_\mu,k_\lambda\> \cr
&=
t^2k(\lambda,\mu)-\phi(\lambda)\bar{\phi}(\mu)\< 
(1-P_{E^\perp })k_\mu, k_\lambda \> \cr                               
&=(t^2-\phi(\lambda)\bar{\phi}(\mu))k(\lambda,\mu) 
+ \phi(\lambda)\bar{\phi}(\mu) \<P_{E^\perp} k_\mu ,\,k_\lambda\>.
 \cr }
$$
That is,
$$
(t^2-\phi(\lambda)\bar{\phi}(\mu))k(\lambda,\mu)=\<(t^2-\MM^*)k_\mu
,\,k_\lambda \> - \phi(\lambda)\bar{\phi}(\mu) \<P_{E^\perp} k_\mu ,\,
k_\lambda \>.
$$
Since $P_{E^\perp}$ has rank $m$, the kernel
$-\phi(\lambda)\bar{\phi}(\mu)
\<P_{E^\perp} k_\mu, k_\lambda\>$ has at most $m$ negative squares. If
$t\geq\|M
\|$
the kernel $\<(t^2-\MM^*)k_\mu,\, k_\lambda\>$ is positive. Thus
$$(t^2-\phi(\lambda)\bar{\phi}(\mu))k(\lambda,\mu)$$ has at most $m$
negative squares as long as $t\geq\|M\|$, in particular when $t =\|M\| =
s_m(\phi)$.

(2)$\implies$ (3). We begin with a simple observation that
explains the introduction of the class $\HH$.

\begin{lemma}
\label{l1.4}
Let $M\in {\cal L}(H)$, let $\phi:D_\phi\subset \Omega \to \C$ be a function
whose domain $D_\phi$ is a set of uniqueness of both $H$ and $\HH$,
and let $F$ be an m-dimensional subspace of $H$, where
$0\leq m < \infty$. Suppose that, for all $h\in F^\perp$ and $\lambda
\in D_\phi $,
$$
(Mh)(\lambda)=\phi(\lambda)h(\lambda).
$$
Then the subset $\cal{S} = \{\lambda:k_\lambda \in F\}$ of $\Omega$
contains at most $m$ points, and there is a function $\psi:\Omega
\setminus\cal{S}\to \C$ such that $\psi(\lambda) =\phi(\lambda)$ for
$\lambda\in D_\phi\setminus \cal{S}$ and 
\begin{equation}
\label{defpsi}
(Mh)(\lambda) = \psi(\lambda)h(\lambda)
\end{equation}
for all $\lambda \in \Omega \setminus \cal{S}$ and $ h\in F^\perp$.
\end{lemma}

\begin{proof}
Since $k$ is nonsingular, the $k_\lambda$ are linearly independent, and
so $\cal{S}$ can contain at most $m$ points. Consider any $\lambda\in 
\Omega\setminus\cal{S}$. Pick an $h\in F^\perp $ such that
$h(\lambda)\neq 0$ (this is possible, since otherwise $h\in F^\perp $
implies $h\in k^\perp_{\lambda}$, i.e., $k_\lambda \in F$). Define
$\psi(\lambda)$ to be $(Mh)(\lambda)/ h(\lambda)$. This
definition is independent of the choice of $h\in F^\perp $ such that
$h(\lambda)\neq 0$. For suppose $h_1, h_2$ are two choices. Then for
any $\mu\in D_\phi$ we have
$$
(h_1\cdot Mh_2)(\mu)-(Mh_1\cdot h_2)(\mu) = h_1(\mu)\phi(\mu)h_2(\mu) -
\phi(\mu) h_1(\mu) h_2(\mu) = 0.
$$
Since $D_\phi$ is a set of uniqueness for $\HH$ it follows that
$h_1 \cdot Mh_2 - Mh_1 \cdot h_2 = 0$, and so
$$
\frac{Mh_1(\lambda)}{h_1(\lambda)} =
\frac{Mh_2(\lambda)}{h_2(\lambda)}\,.
$$
Thus $\psi:\Omega \setminus \cal{S}\to \C$ is well defined.

Now consider $\lambda \in \Omega \setminus \cal{S}$ and 
$g\in F^\perp$. Choose $h\in F^\perp$ such that $h(\lambda)\neq 0$.
The same argument as above shows that $h \cdot Mg-Mh \cdot g = 0$.
Hence
$$
(Mg)(\lambda)=\frac{(Mh)(\lambda)}{h(\lambda)}g(\lambda)=\psi(\lambda)g
(\lambda).
$$
Thus \eqref{defpsi} is satisfied. Pick $\lambda \in D_\phi
\setminus\cal{S}$ and $h\in F^\perp $ such that $h(\lambda)\neq 0$.
Then
$$
\phi (\lambda)h(\lambda ) = (Mh)(\lambda) =
\psi(\lambda)h(\lambda),\text{
and so } \phi(\lambda)= \psi(\lambda).
\eqno\qed
$$
\let\qed\relax
\end{proof}

Example \ref{e1.3a} shows that $\cal{S}$ \emph{can} contain points of
$D_\phi$.

For the remainder of this section denote by $\cal{F}$ the set of
finite subsets of $D_\phi $. For $F\in\cal{F}$, let
$$
\M_F = \span \{k_\lambda:\lambda\in F\}
$$
and let $\cal{R}_\phi $ be the closed linear span of $\{\bar{\phi}
(\lambda)k_\lambda: \lambda\in D_\phi \}$.

Define $T_F\in\cal{L}(\M_F)$ by 
$$T_Fk_\lambda =
\bar{\phi}(\lambda)k_\lambda, \quad\hbox{for }\lambda\in F,
$$ and note the following
standard calculation. For any $x=\sum_{\lambda\in F}c_\lambda
k_\lambda \in \M_F $ we have
$$
\|x\|^2 - \|T_Fx\|^2 =
c^*\bigl((1-\phi(\lambda)\bar{\phi}(\mu))k(\lambda,\mu)\bigr) c, 
$$
where $c=[c_\lambda]_{\lambda\in F}$, regarded as a column vector. It
follows that, if $\phi$ is an $m$-pseudomultiplier and $s_m(\phi) \leq
1$, 
then $1-T_F^*T_F$ has at
most $m$ negative eigenvalues.

\begin{lemma}
\label{l1.5}
Let $\phi$ satisfy condition {\rm 2} in Theorem \ref{textn} and let
$t_0$ be the least
value of $t \geq 0$ such that the kernel \eqref{pickt} has at most $m$
negative squares on $D_\phi$.
There exist functions 
$$
a_j: D_\phi \to \C, \quad\hbox{for }1 \leq j\leq m,
$$
and a positive operator $B$ on $H$ such that $\|B\| = t_0^2$,
$BH \subset {\cal R}_\phi $, and
\begin{equation}
\label{defB}
\<Bk_\lambda, k_\mu \> = \bar{\phi} (\lambda)\phi (\mu)\<
k_\lambda,k_\mu \> - \sum^m_{j=1}\bar{a}_j(\lambda)a_j(\mu) \text{ for
all } \lambda, \mu\in D_\phi.
\end{equation}
\end{lemma}
\begin{proof}
We can assume $t_0=1$. For any $F\in\cal{F}$, the operator
$1-T_F^*T_F$ has at
most $m$ negative eigenvalues. Hence there exists a contraction $B_F$
on $\M_F$ with $0\leq B_F \leq 1$ and vectors $u^F_{1},\ldots, u^F_m \in
\M_F$ such that
\begin{equation}
\label{defBF}
T_F^*T_F = B_F + \sum^m_{j=1}u^F_j \otimes u^F_j.
\end{equation}
The $u^F_j$ can be taken to be suitably scaled Schmidt vectors of
$T_F$ (some could be zero). Note that \eqref{defBF} is equivalent to
\begin{equation}
\label{BF2}
\!\! \bar{\phi}(\lambda)\phi (\mu) \< k_\lambda, k_\mu\> = \< B_F
k_\lambda, k_\mu\> + \sum^m_{j=1} u^F_j(\lambda)^-\, u^F_j(\mu)
\!\!
\end{equation}
for all $F\in\cal{F}$ and $\lambda, \mu \in F $.
In particular,
\begin{equation}
\label{uFbound}
\sum^m_{j=1}|u^F_j(\lambda)|^2 = |\phi(\lambda)|^2\|k_\lambda\|^2 -
\<B_Fk_\lambda ,k_\lambda \>
 \leq |\phi(\lambda)|^2 \|k_\lambda\|^2
\quad\hbox{for }\lambda\in F\in\cal{F}.
\end{equation}
If $\phi(\lambda) = 0$ for some $\lambda$, then
$$
\sum^m_{j=1}|u^F_j (\lambda)|^2 = - \<B_F k_\lambda, k_\lambda\>,
$$
and we must have $B_F k_\lambda =0$. Hence $B_F \M_F \subset
\cal{R}_\phi $ for all $F\in \cal{F}$.

Let $B^\sharp_F = B_F \oplus O_{{\cal M}^\perp_F} \in \cal{L}(H)$. Then
$0\leq B^\sharp_F \leq 1$ and $B^\sharp_F H\subset{\cal R}_\phi $ for
all
$F\in\cal{F}$. For each $F\in\cal{F}$ let
$$
Q_F = \{G\in\cal{F}: F\subset G\}.
$$
No $Q_F$ is empty, and $Q_{F_1} \cap Q_{F_2} = Q_{F_1\cup F_2}$. Thus
$\{Q_F: F\in\cal{F}\}$ is a filter base on $\cal{F}$, generating a
filter
$\cal{Q}$. Let $\cal{U}$ be an ultrafilter on $\cal{F}$ that refines
$\cal{Q}$. Let $\alpha$ be the mapping from $\cal{F}$ to the closed unit
ball $\cal{B}$ of $\cal{L}(H)$ given by $\alpha(F) = B^\sharp_F$.
One can think of $\alpha $ as a continuous mapping from $\cal{F}$ with
the discrete topology to $\cal{B}$ with the weak operator topology.
$\cal{B}$ is compact and so $\alpha$ extends in a unique way to a
continuous mapping from the Stone-\v{C}ech compactification
$\Beta \cal{F}$ of $\cal{F}$ into $\cal{B}$. Since $\cal{U}$ can be
identified with a point of $\Beta \cal{F}$ it follows that
\begin{equation}
\label{gotB}
B=\lim_{F\to \cal{U}} \alpha (F) = \lim_{F\to \cal{U}} B^\sharp_F
\end{equation}
exists in $\cal{B}$. Clearly $0\leq B\leq 1$ and
$BH\subset\cal{R}_\phi$. Recall that $u^F_j(\lambda)$ is defined 
when $\lambda \in D_\phi $, $1\leq j \leq m$ and 
$\lambda \in F \in \cal{F}$. With the same range of $\lambda$ and $j$,
extend the definition to all $F\in \cal{F}$ by setting
$u^F_j(\lambda) = 0$ if $\lambda \not\in F$. For fixed $j$ and
$\lambda \in D_\phi $ the complex-valued function $F\mapsto
u^F_j(\lambda)$ is bounded on $\cal{F}$, by \eqref{uFbound}. Hence
$$
a_j(\lambda)\eqdef \lim_{F\to \cal{U}} u^F_j(\lambda )
$$
exists for all $\lambda \in D_\phi$ and $1\leq j\leq m$. Again by
\eqref{uFbound} we have, for all $\lambda \in D_\phi $,
\begin{equation}
\sum^m_{j=1} |a_j(\lambda)|^2 \leq |\phi (\lambda)|^2 \|k_\lambda\|^2.
\label{abound}
\end{equation}
We now have to show that we can take limits in \eqref{BF2} to obtain
\eqref{defB}.

Fix $\lambda,\mu \in D_\phi $ and let $\epsilon > 0$. By virtue of
\eqref{abound} and \eqref{uFbound} there exists $U\in \cal{U}$ such that
\begin{equation}
\label{diff1}
\left | \sum_j a_j(\lambda)^-\, a_j(\mu) - \sum_j u^F_j(\lambda)^-\,
u^F_j(\mu) \right | < \frac{\epsilon}{2}
\end{equation}
whenever $F\in U$. Similarly, there exists $V\in \cal{U}$ such
that
$$
\bigl | \< Bk_\lambda - B^\sharp_F k_\lambda ,\, k_m \> \bigr | <
\frac{\epsilon}{2}
$$
whenever $F\in V$. Pick $F\in U\cap V\cap Q_{\{\lambda\}}$ (this set
belongs to $\cal{U}$, hence is nonempty). Then $\lambda\in F$, and so
$B^\sharp_F k_\lambda = B_F k_\lambda$. Hence
$$
\left | \< Bk_\lambda - B_F k_\lambda ,\, k_\mu \> \right | <
\frac{\epsilon}{2}.
$$
Since $F\in U$, inequality \eqref{diff1} is also satisfied. Combining
these
inequalities with \eqref{BF2} we obtain
$$
\Bigl | \bar{\phi}(\lambda) \phi (\mu) \< k_\lambda,\, k_\mu \> - 
\< Bk_\lambda,\, k_\mu\> - \sum_j \bar{a}_j(\lambda) a_j(\mu)
\Bigr | < \epsilon.
$$
Since $\epsilon$ was arbitrary, \eqref{defB} is satisfied.

It remains to show that $\|B\| \geq t_0^2$.  By \eqref{defB}, for any $t
\in  \R
$
and $\lambda,~ \mu \in D_\phi$,
$$
(t^2 - \phi(\mu)\bar{\phi}(\lambda))k(\mu,\lambda) =\< (t^2
-B)k_\lambda, k_\mu 
\> -
 \sum^m_{j=1}a_j(\mu)\bar{a}_j(\lambda).
$$
If $t^2 = \|B\|$ then the first term on the right hand side is positive
and the second
has $m$ negative squares.  Thus $(t^2 -
\phi(\mu)\bar{\phi}(\lambda))k(\mu,\lambda)$ has
$m$  negative squares, and so $t_0^2 \leq t^2 = \|B\|.$
\end{proof}

In Example \ref{e1.2a} ($\lambda^{-m}$ on $\D$) 
one can take each $B_F$ to be $1_{\M_F}$, hence
$B=1_H$, and one finds that $a_1(\lambda) = \lambda ^{-m}$. In Example
\ref{e1.3a} (the characteristic function of a point) 
it transpires that $B=0$ and $a_1=\phi $.  In both cases  $a_j$ fails to
belong to $H$. 

\begin{lemma}
\label{l1.6}
For $\phi$ as in Lemma {\rm \ref{l1.5}} there exist $A$ and $L$ in
$\cal{L}(H)$ such that
\begin{enumerate}
\item[(1)] the range and cokernel of both $A$ and $L$ are
contained in $\cal{R}_\phi$;
\item[(2)] $0\leq A\leq 1$ and $\rank A \leq m$;
\item[(3)] $\| L\| = t_0$ and, for every $\lambda\in
D_\phi$ and $h\in H$,  
\begin{equation}
(Lh)(\lambda) = \phi(\lambda)\bigl(h(\lambda) - (Ah)(\lambda) \bigr)
\label{propL}
\end{equation}
and
\begin{equation}
\label{propL*}
L^*k_\lambda = (1-A)\bar{\phi}(\lambda) k_\lambda.
\end{equation}
\end{enumerate}
\end{lemma}
\begin{proof}
The relation \eqref{defB} in Lemma \ref{l1.5} can be written
$$
\< \bar{\phi} (\lambda)k_\lambda,\, \bar{\phi}(\mu)k_\mu \>_H =
\<B^{{1}/{2}}k_\lambda \oplus (\bar{a}_1(\lambda), \dots,
\bar{a}_m(\lambda)),\,
 B^{{1}/{2}}k_\mu \oplus (\bar{a}_1(\mu)\dots, \bar{a}_m(\mu)) \> _{H
\oplus \C^m}
$$
for all $\lambda, \mu \in D_\phi$. Hence we can define an isometry
$V: {\cal R}_\phi \to {\cal R}_\phi \oplus \C^m$ by 
\begin{equation}
\label{defV}
V\bar{\phi} (\lambda)k_\lambda = B^{{1}/{2}}k_\lambda \oplus
(\bar{a}_1(\lambda), \ldots,\bar{a}_m(\lambda))\quad
\hbox{for } \lambda \in D_\phi.
\end{equation}
The isometry $V$ necessarily has the form
\begin{equation}
\label{formV}
V = \left [
\begin{matrix}
C\\
1\otimes f_1\\
\vdots\,\\
1\otimes f_{\rlap{$\scriptstyle m$}\phantom1}\\
\end{matrix}
\right ]
\end{equation}
for some $C: {\cal R}_\phi \to {\cal R}_\phi $ and $f_1, \ldots, f_m \in
{\cal R
}_\phi$.
The relation $V^* V=1$ yields
$$
C^* C + f_1\otimes f_1 + \cdots + f_m \otimes f_m = 1.
$$
Let
$$
A= f_1 \otimes f_1 + \cdots + f_m \otimes f_m.
$$
Then $0 \leq A \leq 1$, $AH\subset{\cal R}_\phi$, $\rank A \leq m$ and
$C^* C
= 1-A$. From \eqref{defV} and \eqref{formV} we have
\begin{equation}
\label{defC}
\eqaligntop{
C \bar{\phi} (\lambda)k_\lambda &= B^{1/2} k_\lambda\cr
\bar{\phi}(\lambda)\bar{f}_j(\lambda) &= \bar{a}_j(\lambda)}
\end{equation}
for $\lambda \in D_\phi$ and $1 \leq j \leq m$. Extend $C$ to an
operator
$C^\sharp$ on $H$ by taking the direct sum of $C$ and the zero
operator on ${\cal R}^\perp_\phi$, and define
\begin{equation}
\label{defL}
L = B^{1/2}C^\sharp : H\to H.
\end{equation}
Clearly the range and cokernel of $L$ are contained in ${\cal R}_\phi$.
Since $\|B^{1/2}\| = t_0$ and $\|C\| \leq 1$ we have
$$
\|L\| \leq t_0.
$$
On applying $C^*$ to \eqref{defC} we find that, for $\lambda\in D_\phi$,
$$
L^*k_\lambda = (C^\sharp)^* B^{1/2}k_\lambda = C^*
B^{1/2}k_\lambda = C^*C \bar{\phi} (\lambda)k_\lambda =
(1-A)\bar{\phi}(\lambda) k_\lambda .
$$
Hence, for any $h\in H$ and $\lambda \in D_\phi$, 
$$
\eqalignbot {
(Lh)(\lambda) &=  \<Lh, k_\lambda\> = \<h,L^*k_\lambda\> 
               =   \<h, (1-A) \bar{\phi}(\lambda)k_\lambda\> \\
              &=  \phi(\lambda) \<(1-A)h,\, k_\lambda\> 
               =  \phi(\lambda)\bigl(h(\lambda)-(Ah)(\lambda)\bigr). \\ 
}
$$
It remains to show that $\|L\| \geq t_0,$ or equivalently, that the
kernel
$$
\left(\|L\|^2 - \bar{\phi}(\lambda) \phi(\mu)\right) k(\mu, \lambda)
$$
has at most $m$ negative squares.  From \eqref{propL*} we have, for
$\lambda, \mu \in
\D_\phi,$
$$
\eqalignbot {
\< L^*k_\lambda, L^*k_\mu \> &= \<(1-A)\bar{\phi}(\lambda)k_\lambda,
                        (1-A)\bar{\phi}(\mu)k_\mu \>  \\
                        &= \bar{\phi}(\lambda) \phi(\mu) \<
(1-2A+A^2)k_\lambda,
k_\mu \>.
}
$$
Hence
$$
\eqalignbot {
\left(\|L\|^2 - \bar{\phi}(\lambda) \phi(\mu)\right) k(\mu, \lambda)
        &= \< (\|L\|^2 -LL^*) k_\lambda, k_\mu \> + \\
        & \qquad \bar{\phi}(\lambda) \phi(\mu)\left[ \<
(1-2A+A^2)k_\lambda,k_\mu \>
                - \<k_\lambda, k_\mu \> \right]  \\
        &= \text{a positive kernel} - \bar{\phi}(\lambda) \phi(\mu)
                 \< A(2-A)k_\lambda,k_\mu \>.
}
$$
Since $A(2-A)$ is positive and has rank $\leq m$, the right hand side
has at most
$m$ negative squares, as required.
\end{proof}

We can now prove that (2) $\Rightarrow$ (3) in Theorem \ref{textn}. 
Suppose that $\phi$ satisfies (2), the least value of $t$ being $t_0$. 
By Lemma \ref{l1.6} there exist $A$ and $L$ in $L(H)$ such that
$||L||=t_0$ and, for all $\lambda \in D_{\phi}$ and $h
\in H$,
$$
Lh(\lambda) = \phi(\lambda) (h(\lambda) -(Ah)(\lambda)).
$$
Moreover $A$ has rank $m$, so that {\rm Ker} $A$ has codimension $m$, and
for
all $h \in $ {\rm Ker} $A$ and $\lambda \in D_{\phi}$,
$$
Lh(\lambda) = \phi(\lambda)h(\lambda).
$$
By Lemma \ref{l1.4} there is a subset $\cal{S}$ of $\Omega$ containing
at
most $m$ points and a function $\psi: \Omega \setminus \cal{S}
\rightarrow \C$
such
that $\psi(\lambda)= \phi(\lambda)$ for $\lambda \in D_{\phi}
\setminus \cal{S}$ and
\begin{equation}
 \label{psiL} 
\psi(\lambda)h(\lambda) = Lh(\lambda) =
\phi(\lambda)h(\lambda)
\end{equation}
for all $h \in$ {\rm Ker} $A$ and $\lambda \in \Omega \setminus
\cal{S}$.
  Extend $\psi$ to a function on $D_{\psi}=(\Omega \setminus \cal{S}) \cup
D_{\phi}$ by defining $\psi(\lambda) = \phi(\lambda)$ for $\lambda
\in D_{\phi} \cap \cal{S}$.  Then $\psi$ extends $\phi$ and its
domain $D_{\psi}$ is a set of uniqueness for $H$ (since $D_\psi \supset
D_\phi$).  There are at most $m$ points in $\Omega \setminus D_{\psi}
  \subset \cal{S}$,
and (\ref{psiL}) tells us that $\psi${\rm Ker} $A \subset H$.  
Thus $\psi$ is an
$m$-\ps/~of $H$ which extends $\phi$.  Hence (2) $\Rightarrow$ (3).

Now suppose (1)-(3) hold.  From (\ref{psiL}) we have
$$
s_m(\psi) \le ||L|| = t_0.
$$
In the proof of (1) $\Rightarrow$ (2) we showed that the kernel
\eqref{pickt}
has at most $m$ negative squares if $t=s_m(\phi)$, that is, 
$$
t_0 \leq s_m(\phi).
$$
If $E$ is an $m$-codimensional subspace 
of $H$ such that $\psi E \subset
H$ then $\phi E|D_{\phi} \subset H|D_{\phi}$ 
and $||M_{\psi,E} || =
|| M_{\phi, E}|D_{\phi} ||$.  Taking infima over closed
$m$-codimensional subspaces $E$ such that $\psi E \subset H$ we obtain
$$
s_m(\phi) \le s_m(\psi).
$$
Combining the last three inequalities, we have
$$
s_m(\phi) = t_0 = s_m(\psi).
$$
\hfill $\qed$

We can now derive a simple but informative multiplicity result
for \ps/s.
\begin{theorem}
\label{mplcty}
Let $H$ be a Hilbert space of 
functions with nonsingular kernel on a set $\O$ and let
$\phi$ be an $m$-\ps/ of $H$.  For any $\xi \in \C$  such that
$| \xi | > s_m(\phi)$, the equation $\phi(\lambda) = \xi$
has at most $m$ solutions for $\lambda \in D_\phi$.
\end{theorem}
\begin{proof}
Pick $L,~A$ as in Lemma \ref{l1.6}.  Suppose there are $n$
points $\lambda_1,\dots,\lambda_n \in D_\phi$ which are
solutions of the equation, where $n > m$.  Since $A$ has rank
at most $m$ there exist scalars $c_1,\dots,c_n$, not all zero, such that
$\sum c_j Ak_{\lambda_j} = 0$.
Let $h=\sum c_j k_{\lambda_j}$.  Then $h \neq 0$, and
$$
L^*h = \sum c_j L^* k_{\lambda_j} = \sum c_j \bar{\xi}(k_{\lambda_j}
-Ak_{\lambda_j}) =\bar{\xi}h,
$$
and so $\bar{\xi}$ is an eigenvalue of $L^*$.  However
$$
||L^*|| = s_m(\phi) < |\xi|,
$$
a contradiction.
\end{proof}

We conclude this section with an aside about the key Lemma \ref{l1.6}.
How should we understand the rank $m$ Hermitian operator $A$ associated with
an $m$-\ps/ $\phi$?  Let us consider the case $m=1$, and put $A=f \otimes f$
for some $f \in H$.  The relation \eqref{propL} becomes
\begin{equation}
\label{propL1}
(Lh)(\lambda) = \phi(\lambda)\left(h(\lambda) - \< h,f \> f(\lambda) \right)
\end{equation}
for $h \in H,~\lambda \in D_\phi$.  What is the connection between the \ps/
$\phi$ and the function $f$?  We propose that $f$ be regarded as a singular
vector of the (unbounded) operator $M_\phi$ on $H$ corresponding to the
singular value $s_0(\phi) =\infty$.  Strictly speaking 
this does not make sense, but we can give it a meaning as follows.

Consider first the case of a finite kernel $k$ on $\O = \{ \lambda_1, \dots,
\lambda_n \}$, and define $T$ on $H$ by 
$Tk_\lambda = \bar{\phi}(\lambda)k_\lambda$ (so that $T=T_\O$ in the
notation above).  Then $T$ has a maximising vector $u$ such that, for some
contraction $B$ on $H$,
$$
T^*T = B + u \otimes u.
$$
Here $|| u || = || T ||$, and in the notation of Lemma \ref{l1.5},
$u = a_1$.  By equation \eqref{defC}, $\phi f = u$, and so
$$
f = \frac{1}{\phi} u = (T^*)^{-1} u.
$$
Thus $u = T^* f$.  Hence $f$ is a maximising vector of $T^*$
and
$|| f || = 1$. In other words $f$ is a singular vector of the multiplication
operator $M_\phi$ corresponding to $s_0(\phi)$.  For finite kernels, then,
our statement about $f$ is meaningful and true.

Now consider a nonsingular kernel on a general set $\O$ and apply the
foregoing observation to $T_F$ for each $F \in \cal F$.  Then there is a unit
maximising vector $f_F \in {\cal M}_F$ and a contraction $B_F $ on ${\cal M}_F$ 
such that 
$$
T_F^* T_F = B_F + T_F^*f_F \otimes  T_F^*f_F.
$$
Along the ultrafilter $\cal U$, $f_F$ tends weakly to a limit,
which must equal $f$.  If $|| f || < 1$ then $\phi$ is a multiplier
of $H$ (put $h=f$ in \eqref{propL1}).  If  $|| f || = 1$ then we have
$f_F \to f$ in norm as $F \to \cal U$.  That is, $f$ is the norm limit of
unit maximising vectors of the approximating finite rank compressions
$T_F^*$, $F \in \cal U$ of the multiplication operator $M_\phi$.

\section{Pseudomultipliers of analytic kernels}

What are the pseudomultipliers of the familiar Hilbert spaces of
analytic functions? The examples we began with ($1/z$ and the
characteristic function of a point) might lead one to the optimistic
hope that the $1$-pseudomultipliers are simply 
the multipliers modified by a
removable singularity or a simple pole. We already know that, for a
general analytic kernel, this can fail to be true in at least three
ways: the domain $\Omega $ could be too small (Example
\ref{ewrongO}), $H$ could really be a space of functions of $z^2$
(Example \ref{e1.7}) or there could be more subtle structural reasons
(Example \ref{eHRS}). Nevertheless, it \emph{is} true for a wide
class of kernels, as we show in Theorem \ref{t2.2} below.

Let $\Omega $ be a domain in $\C$ and let $H$ be a Hilbert space of
analytic functions on $\Omega$ with reproducing kernel $k$.
Such a kernel will be called an {\em  analytic
kernel}.  There is one simple observation we can make about
\ps/s of such kernels.
\begin{theorem}
\label{anker}
Let $k$ be a nonsingular analytic kernel on a domain $\O$ and let $\phi$
be an $m$-\ps/ of the associated Hilbert space $H$ of functions on $\O$.
Let the defect of $\phi$ be $d$.  There is a set
$\cal S \subset D_\phi$, containing at most $m-d$ points,
such that $\phi$ is analytic at every point of $D_\phi \setminus \cal S$.
Moreover, each point of $\O \setminus D_\phi$ is a pole of $\phi$
and each point of $\cal S$ is either a removable singularity or
a pole of $\phi$.
\end{theorem}
Note that it {\em can} happen that a point of $D_\phi$ is a pole of
$\phi$: recall  Example \ref{e1.2b}.
\begin{proof}
By Lemma \ref{l1.6} there exist bounded linear operators $A,~L$ on $H$
such that $0 \leq A \leq 1$, $A$ has rank at most $m$ and, 
for every $h \in H$ and $ \lambda \in D_\phi$,
$$
(Lh)(\lambda) = \phi(\lambda)\left(h(\lambda) - (Ah)(\lambda)\right).
$$
We claim that $\phi$ is analytic at any $\beta \in \O$ such that $k_\beta
\not \in AH$.  Indeed, for such $\beta$,  ${\rm Ker}~ A$ is not a subset of 
$k_\beta^\perp$, and so there exists $h \in H$ such that $Ah =0$ and
$h(\beta) \neq 0$.  We have $Lh = \phi h$, and so $\phi = (Lh)/h$ on
any neighbourhood of $\beta$ on which $h$ does not vanish.  Thus $\phi $
is analytic at $\beta$.

It follows that $k_\beta \in AH$ for the $d$ points of $\O \setminus D_\phi$.
Since $k$ is nonsingular and $AH$ has dimension $\leq m$, the set
$$
{\cal S} \stackrel{\rm def}{=} \{ \beta \in D_\phi: k_\beta \in AH \}
$$
can contain at most $m-d$ points.  The preceding paragraph shows that
$\phi$ is analytic at each point of $ D_\phi \setminus \cal S$.

Now consider any point $\alpha \in\O \setminus D_\phi $. 
Then $\alpha$ is a singularity of $\phi$.  That is, there
exist $u,~g \in H$ such that $u(\alpha) =0,~ \phi u = g|D_\phi$ and
$g(\alpha) \neq 0$.  Clearly $u \neq 0$, so that $\alpha$ is an isolated
zero of $u$.  Hence we can write $\phi = g/u$ on a punctured
neighbourhood of $\alpha$, and so $\phi$ has a pole at $\alpha$.

Consider any point $\alpha \in \cal S$.  Since $\phi$ is analytic in a
punctured neighbourhood of $\alpha$ it must be the case that $\alpha $ is
either a removable singularity, a pole or an essential singularity of
$\phi$.  In the latter case, by the great Picard theorem, $\phi$ attains
all but one complex values infinitely often, contradicting
Theorem \ref{mplcty}.  Hence one of the other two alternatives holds.

\end {proof}

This is an opportune moment to point out that \ps/s of Hilbert spaces of
analytic functions of several complex variables are not interesting objects,
at least as we have defined them here.  By exactly the argument above, an
$m$-\ps/ $\phi$ of such a space is analytic at all but at most $m$ points.
Now a singularity of $\phi$ in the sense of Definition \ref{d1.1} would
have to be an isolated pole of $\phi$, and an analytic function of
several variables cannot have any such.  Thus the $m$ points where analyticity
fails are removable singularities, and $\phi$ is just an analytic function
plus some point discontinuities.  However, it is possible the notion 
of \ps/ could be developed for Hilbert spaces of {\em vector-valued}
functions so as to apply to analytic functions of several variables.

\begin{corollary}
\label{centire}
Let $H$ be a Hilbert space of entire functions having nonsingular kernel 
on $\C$, and let $\phi$ be an $m$-\ps/ of defect $d$ on $H$. 
There exists a rational function $\psi$ of degree at most $m$ which
agrees with $\phi$ at all except at most $m-d$ points of $D_\phi$.
\end{corollary}
Here the degree of $\psi$ is defined to be the sum of the multiplicities 
of the poles of $\psi$, including $\infty$ if applicable.  The points at
which $\phi$ and $\psi$ differ can include poles of $\psi$.
\begin{proof}
By Theorem \ref{anker} there is a set ${\cal S} \subset D_\phi$
containing at most $m-d$ points, consisting of poles and removable
singularities of $\phi$, such that $\phi$ is analytic on $D_\phi \setminus
\cal S$.  Let $\psi$ be the function obtained from $\phi$ by removing any
removable singularities and deleting from the domain of $\phi$ any point
$\alpha \in \cal S$ which is a pole of $\phi|D_\phi \setminus \{\alpha\}$.
Then $\psi$ is meromorphic in the whole complex plane.  Again 
 by the great Picard theorem and
Theorem \ref{mplcty}, $\psi$ cannot have an essential singularity at $\infty$
and hence is rational.  Clearly $\psi$ agrees with $\phi$ except at points of
$\cal S$.
\end{proof}

In order to get more detailed information 
about the nature of \ps/s of analytic kernels we
 consider kernels satisfying the following two axioms:
\begin{description}
\item[\hskip-\parindent\normalfont (AK1)] $H$ is invariant under the
operation $M_z$ of multiplication by
the independent variable;
\item[\hskip-\parindent\normalfont (AK2)] Every bounded linear operator
on $H$ which commutes with $M_z$ is a multiplication operator.
\end{description}

Not all spaces of interest have these properties (e.g., the Fock space
does not satisfy (AK1)), but many of them do. An easy consequence of
(AK1) is
\begin{equation}
\label{eigen}
M^*_z k_\lambda = \bar{\lambda}k_\lambda 
\end{equation} 
for all $\lambda \in \Omega$.
 
{\bf Remark.}
{\normalfont (AK2)} 
holds if {\normalfont (AK1)} is true and 
either $H$ contains the polynomials on $\O$ as a dense subspace
or every
 eigenvector of $M^*_z$
has the form $ck_\alpha$ for some
$c\in \C$ and $\alpha\in \Omega $.

For suppose $TM_z = M_zT$.
In the former case it is plain that $T$ is multiplication by $M 
\boldsymbol{1}$.  In the latter we have,
 for any $\lambda \in \Omega $,
$$
M^*_z T^*k_\lambda = T^* M^*_zk_\lambda = \bar{\lambda}T^* k_\lambda 
$$
and so $T^*k_\lambda$ is either $0$ or an eigenvector of $M^*_z$ with
eigenvalue $\bar{\lambda}$. Hence, by assumption,
$$
T^* k_\lambda = \bar{\psi}(\lambda) k_\lambda
$$
for some $\psi(\lambda) \in \C$. For any $h\in H$ and $\lambda \in
\Omega$
$$
(Th)(\lambda) = \< Th, k_\lambda \> = \<h, T^* k_\lambda\>
                  = \< h, \bar{\psi}(\lambda)k_\lambda\> =
\psi(\lambda)h(\lambda) .
$$
Thus $\psi$ is a multiplier and $T=M_\psi$.

\begin{theorem} 
\label{t2.2}
Let $H$ be a Hilbert space of analytic functions on a domain $\O$.
Suppose that $H$ has nonsingular kernel and that
 {\normalfont (AK1)} and {\normalfont (AK2)}  hold. Let
 $\phi : D_\phi \subset \Omega \to
\C$ be a $1$-pseudomultiplier of $H$. Then there exists $\alpha\in
\C$ and a multiplier $\theta$ of $H$ such that
 $D_\phi$ is either $\O$ or $\Omega \setminus \{\alpha\}$ and
$$
\phi (\lambda) = \frac{\theta (\lambda )}{\lambda - \alpha} \quad\text {
for }
\lambda\in \Omega \setminus \{\alpha\}.
$$ 
\end{theorem}

We do not assert that $\alpha \in \O$ in general.
Observe that the conclusion holds when $\phi$ is a multiplier or even
a function that differs from a multiplier at a single point $\alpha \in \O$
 (take $\theta(\lambda) = (\lambda - \alpha)\phi(\lambda)$).

\begin{proof} 
By Lemma \ref{l1.6} there exist $L \in \cal{L}(H)$ such that $\|L\| =
s_1(\phi)$ and $f\in H$ such that $\|f\| \leq 1$ and, for all $\lambda
\in D_\phi$,
$$
L^*k_\lambda = \bar{\phi}(\lambda )k_\lambda -
\bar{\phi}(\lambda)\bar{f}(\lambda)f
$$
(we have put $A = f \otimes f$).
From this equation and (\ref{eigen}) it follows that 
\begin{equation}
\label{commut}
(M^*_z L^* - L^* M^*_z) k_\lambda =
\bar{\phi}(\lambda)\bar{f}(\lambda)(\bar{\lambda}f - M^*_z f)
\end{equation}
for all $\lambda\in D_\phi$. Hence
$$
\range (M^*_z L^* - L^* M^*_z) \subset \span\{f, M^*_z f\}.
$$
We consider the three possible dimensions of this span. 

Suppose $f=0$: then by
\ref{propL}, for
$\lambda\in D_\phi$ and $h\in H$,
$$
 (Lh)(\lambda) = \phi (\lambda)h(\lambda).
$$
By Lemma \ref{l1.4} and 
Theorem \ref{mltplr} $D_\phi = \O$ and $\phi$ is a multiplier
of $H$.  (Incidentally,
we have in this case
$$
s_0(\phi) = ||M_\phi|| = ||L|| = s_1(\phi),
$$
so that $s_0(\phi)$ is a multiple singular value.)

Secondly, suppose $f$ and $M^*_zf$ span a space of dimension one:
say $M_z^* f = \bar{\alpha} f$ for some $\alpha \in \C$.
In equation \eqref{propL} replace $h$ by $(M_z - \alpha)h$ to
obtain, for all $\lambda\in
D_\phi$ and $h \in H$,
\begin{eqnarray}
(L(M_z - \alpha)h)(\lambda) &=& \phi(\lambda) \left(
(\lambda - \alpha)h(\lambda) - \< (M_z - \alpha)h,f \> f(\lambda)\right)
\nonumber \\
 &=& \phi(\lambda) \left(
(\lambda - \alpha)h(\lambda) - \<h, M_z^*f  - \bar{\alpha} f\>
 f(\lambda)\right)\nonumber \\
&=&(\lambda - \alpha) \phi(\lambda)h(\lambda).\nonumber
\end{eqnarray}
Since $L(M_z - \alpha)$ is bounded,
Lemma \ref{l1.4} tells us that there is a 
multiplier $\theta$ of $H$ such that
\begin{equation}
\label{defth}
\theta(\lambda) = (\lambda-\alpha)\phi(\lambda) \text{ for all }
\lambda \in D_\phi.
\end{equation}
We claim that $\phi $ can have no singularity in $\O$ other than $\alpha$.
For suppose $\beta \neq \alpha$ is a singularity.  Then there exists
$h \in H$ such that $h(\beta) =0$ and $\phi h$ has an extension $u \in H$
such that $u(\beta) \neq 0$.  Now $(M_z - \alpha) u \in H$ is an extension
of the function $\lambda \mapsto (\lambda - \alpha) (\phi h)(\lambda)$ on 
$D_\phi$, that is, of $\theta h |D_\phi$.  Since $\theta h \in H$, it
follows that 
 $\theta h = (M_z - \alpha) u$.  Evaluating both sides at $\beta$ we have
$$
0=\theta(\beta) h(\beta) =(\beta - \alpha)u(\beta) \neq 0,
$$
a contradiction.  Thus $\beta$ cannot be a singularity, and so $D_\phi$
is either $\O$ or $\O \setminus \{\alpha\}$.  It follows that
 $D_\phi \setminus \{\alpha\} =\O \setminus \{\alpha\}$, and, from
\eqref{defth},
$$
\phi (\lambda) = \frac{\theta (\lambda )}{\lambda - \alpha} \quad\text {
for }
\lambda\in \Omega \setminus \{\alpha\}
$$ 
as required.

The remaining possibility is that $\span \{f, M^*_zf\}$ has dimension
2. Then, by (\ref{commut}), there exist $u,v\in H$ such that
\begin{equation}
\label{uv}
M^*_z L^* - L^*M^*_z = f\otimes u - M^*_zf \otimes v.
\end{equation}
On applying both sides to $k_\lambda$ we find, for $\lambda\in
D_\phi$,
$$
\bar{\phi}(\lambda)\bar{f}(\lambda)(\bar{\lambda}f - M^*_z f) =
\bar{u}(\lambda)f - \bar{v}(\lambda) M^*_z f.
$$
By the linear independence of $f$ and $M^*_zf$,
$$
\bar{\lambda}\bar{\phi}(\lambda) \bar{f}(\lambda) = \bar{u} (\lambda),
\quad \bar{\phi}(\lambda)  \bar{f}(\lambda) = \bar{v}(\lambda)
$$
on $D_\phi$. It follows that $u-M_zv$ vanishes on $D_\phi$, hence is
zero. Thus (\ref{uv}) becomes
$$
M^*_z L^* - L^* M^*_z = (f \otimes v) M^*_z - M^*_z(f \otimes v).
$$
Consequently $L+v \otimes f$ commutes with $M_z$. Thus there exists a
multiplier $\theta$ of $H$ such that $L + v \otimes f = M_\theta$.
Then, for $\lambda \in D_\phi$,
$$
\bar{\theta }(\lambda)k_\lambda  = M^*_\theta k_\lambda =
L^*k_\lambda + (f \otimes v) k_\lambda
                              =\bar{\phi}(\lambda)k_\lambda -
\bar{\phi}(\lambda)\bar{f}(\lambda)f + \bar{v}(\lambda)f
                              =\bar{\phi}(\lambda)k_\lambda.
$$
Thus $\phi$ is the restriction of a
multiplier.  By Theorem \ref{mltplr},
$\phi$ {\em is} a multiplier.
\end{proof}

It would be of interest to identify the \ps/s and the 
corresponding quantities  $s_m(\cdot)$ for favourite function spaces.
As we have mentioned, in the case of $H^2$, a well known theorem of
Adamyan, Arov and Krein provides a complete description.  To present it
we introduce some terminology.  Suppose that $\phi,~\psi$ are functions with
domains  $D_\phi,~D_\psi$ contained in a set $\O$.  We shall say that $\psi$
is a {\em finite modification} of $\phi$ if $D_\psi \supset D_\phi,~
D_\psi \setminus D_\phi$ is finite and $\phi,~\psi$ agree at all but finitely
many points of $D_\phi$.  We recall that $H^\infty_{(\ell)}$ denotes the
set of functions of the form $\phi=f/p$ where $f \in H^\infty$ and $p$ is a
polynomial with at most $\ell$ zeros in $\D$, counting multiplicities,
and none of unit modulus; we take
$D_\phi$ to be the complement of the set of poles of $\phi$ in $\D$
(thus we suppose that all removable singularities have been removed).
\begin{theorem}
\label{taak}
A function $\phi: D_\phi \subset \O \to \C$ is a \ps/ of 
 $H^2$ if and only if there is a function $\psi \in
H^\infty_{(\ell)}$, for some non-negative integer $\ell$, such that
$\phi$ is a finite modification of $\psi$.
Moreover, if $\phi$ is an $m$-\ps/ of $H^2$, then $s_m(\phi)=
\| \phi \|_{L^\infty({\mathbb T})}$.
\end{theorem}
This is essentially the main result of [AAK], where however it is expressed
in a rather different context.  A careful derivation of this form
of the theorem from the original one can be found in [Q2, Sec. 7].

We shall conclude with a description of the \ps/s of the Fock space
(see Example \ref{e1.3b}).
Recall that a rational function is said to be {\em proper} if it has a finite
limit at infinity.
\begin{theorem}
\label{tfock}
The \ps/s of the Fock space $\Phi$ are precisely the finite modifications 
of the proper rational functions.
\end{theorem}

\begin{proof} Let $\cal E$ be the space of all entire functions.  We claim that
if $p$ is a non-zero polynomial, $f$ is an entire function and $pf \in \Phi$
then $f \in \Phi$.  Indeed, if $p$ is non-constant, 
$\{z: |p(z)| < 1\}$ is a bounded subset of $\C$, and
so
$$ \eqalign{
\int_\C |f(z)|^2 e^{-|z|^2} dA &= \int_{|p(z)|<1} + \int_{|p(z)|>1} 
|f(z)|^2 e^{-|z|^2} dA \cr
&\leq
\int_{|p(z)|<1}|f(z)|^2 e^{-|z|^2} dA +
\int_{|p(z)|>1} |p(z) f(z)|^2 e^{-|z|^2} dA < \infty.
}
$$

Let $\phi$ be a proper rational function of degree $n$ with denominator $p$.
We show that $\phi$ is a \ps/ of $\Phi$.  Let $E = \Phi \cap p\cal E$.  Then
$E$ has codimension $n$ in $\Phi$: it is given by the vanishing of $n$ linear
functionals, which are evaluations of $f$ and some of its derivatives at
zeros of $p$.  Write $\phi$  as a sum of partial 
fractions in which each term is
either constant or is a constant divided by a divisor of $p$. 
The claim above shows that each such term multiplies $E$ into
$ \Phi$ and hence $\phi E \subset \Phi$. 
 It is clear that each pole of $\phi$ is
a singularity in the sense of Definition \ref{d1.1}, and hence $\phi$ is an
$n$-\ps/, whose domain $D_\phi$ is the complement of the set of poles of $\phi$.
Now consider the function $\psi$ obtained by changing the value of $\phi$ at
points $\alpha_1,\dots,\alpha_r \in D_\phi$ and giving it a value at poles 
$\beta_1,\dots,\beta_s$ of $\phi$.
Then $\psi$ has domain $D_\psi = D_\phi \cup \{\beta_1,\dots,\beta_s\}$.
  If $\beta_j$  has multiplicity $n_j$
then $\psi$ multiplies the closed $(n+r+s)$-codimensional subspace
$$
\{f \in \Phi \cap p{\cal E} : f(\alpha_j)=0,~1 \leq j \leq r,~
f^{(n_i)}(\beta_i)=0,~1 \leq i \leq s \}
$$
of $\Phi$ into $\Phi|D_\psi$.  
Points outside the domain of $\psi$ are poles of
$\phi$ and they remain singularities of $\psi$.  
Thus $\psi$ is an $(n+r+s)$-\ps/ of $\Phi$. That is, any finite modification of
a proper rational function is a \ps/ of $\Phi$.

Conversely, consider an $m$-\ps/ $\psi$ of $\Phi$.  Let $\psi$ have defect $d$. 
By Corollary \ref{centire} there exists a rational function $\phi$ which
differs from $\psi$ at no more than $m-d$ points.  We have to show that 
$\phi$ is proper.

The function $z$ is not a \ps/ of $\Phi$.  For suppose $s_m(z) = t < \infty$.
Then  by Theorem \ref{textn} the kernel $(t^2 - \lambda \bar{\mu})
e^{\lambda \bar{\mu}}$ has at most $m$ negative squares.  However,
$$
(t^2 - w) e^w = t^2 + \sum_{n=1}^\infty \left(\frac{t^2}{n} -1 \right)
\frac{w^n}{(n-1)!},
$$
which has infinitely many negative coefficients, a contradiction.  It follows
that $z \mapsto z-a$ is not a \ps/, for any $a \in \C$.  In fact no non-constant
polynomial is a \ps/ of $\Phi$.  For suppose that $p$ is such a polynomial and
$pE \subset \Phi$ for some closed finite-codimensional space $E \subset \Phi$.
Factorize $p$ as $p(z)= (z-a)q(z)$ for some $a \in \C$ and polynomial $q$.
By the claim above, if $(z-a)qE \subset \Phi$ then $(z-a)E \subset \Phi$,
contradicting the fact that $z-a$ is not a \ps/.    Hence $p$ is not a \ps/.

Return to the rational \ps/ $\phi$ of $\Phi$ and expand it in partial fractions:
$$\phi = p + \sum u_j
$$
where 
$p$ is a polynomial and each $u_j$ is a 
constant divided by a power of a linear
function.  We showed above that each $u_j$ is a \ps/, and so $\phi -\sum u_j$
is a \ps/.  That is, $p$ is a \ps/.  It follows that $p$ is constant, and hence
that $\phi$ is proper.
\end{proof}


\begin{thebibliography}{}
\bibitem[AAK]{AAK}
V. M.  Adamyan, V. Z.  Arov and M. G.  Krein, Analytic properties of
Schmidt pairs of a Hankel operator and generalized Schur--Takagi
problem, {\it Mat.  Sbornik} {\bf 86} (1971), 33--73. 

\bibitem[Ag1]{Ag1}
J. Agler, Interpolation, preprint (1987).

\bibitem[Ag2]{Ag2} 
J. Agler, Nevanlinna--Pick interpolation on Sobolev space,  
{\it Proc. 
Amer. Math. Soc.} {\bf 108}(1990) 341-351.

\bibitem[Ar]{Ar}
N. Aronszajn, Theory of reproducing kernels, {\it Trans. Amer. Math.
Soc.}
{\bf 68} (1950), 337--404.

\bibitem[BB]{BB}
F. Beatrous and J. Burbea, Positive definiteness and its applications to
interpolation
problems for holomorphic functions, {\it Trans. Amer. Math. Soc.},{\bf
284} (1984) 247--270.

\bibitem[H]{H} H. Hedenmalm, An invariant subspace of  the Bergman space
having the
codimension two property, {\it J. Reine Angew. Math.} {\bf 443} (1993)
1--9.

\bibitem[HRS]{HRS} H. Hedenmalm, S. Richter and K. Seip, Interpolating
sequences and
invariant subspaces of given index in the Bergman space, preprint 1995.

\bibitem[P]{P}  
G.  Pick, \"Uber die Beschr\"ankungen analytischer Funktionen, welche
durch vorgegebene Funktionswerte bewirkt werden, {\it Math.  Ann.}  {\bf
77}
(1961), 7--23.

\bibitem[Q1]{Q1}
P. Quiggin, For which reproducing kernel Hilbert spaces is Pick's
theorem true? {\it Integral Equations and Operator Theory\/} {\bf 16}
(1993), 244--266.

\bibitem[Q2]{Q2}
P. Quiggin, Generalisations of Pick's Theorem to Reproducing Kernel
Hilbert Spaces,
 Ph. D. thesis, Lancaster University, 1994.

\bibitem[S]{S}
F. Szafraniec, On bounded holomorphic interpolation in several
variables, {\it
Monatshefte Math.} {\bf 101} (1986) 59-66.

\end{thebibliography}
\end{document}